\documentclass[runningheads]{llncs}
\usepackage{graphicx}

\usepackage{amssymb,amsmath}
\usepackage{xcolor}
\usepackage{hyperref}

\begin{document}
\title{A note on the $p$-operator space structure of the $p$-analog of the Fourier-Stieltjes algebra}
\titlerunning{$p$-operator space structure of $B_p(G)$}
%

\author{Mohammad Ali Ahmadpoor\inst{1}\orcidID{0000-0001-6902-1916} \and Marzieh Shams Yousefi \inst{2}\orcidID{0000-0003-0426-708X}}
\authorrunning{M. A. Ahmadpoor \& M. Shams Yousefi}
%
\institute{${}^{{1,2}}$ Department of Pure Mathematics, Faculty of Mathematical Sciences, University of Guilan, Rasht, Iran\\
\email{m-a-ahmadpoor@phd.guilan.ac.ir}\\
\email{m.shams@guilan.ac.ir}}

\maketitle              
\begin{abstract}
In this paper one of the possible $p$-operator space structures of the $p$-analog of the Fourier-Stieltjes algebra  will be introduced, and to some extend will be studied. This special sort of operator structure will be given from the predual of this Fourier type algebra, that is the algebra of universal $p$-pseudofunctions. Furthermore, some applicable and expected results will be proven. 

Current paper can be considered as a new gate into the collection of problems around the $p$-analog of the Fourier-Stieltjes algebra, in the $p$-operator space structure point of view.

\keywords{$p$-Operator spaces\and $p$-Analog of the Fourier-Stieltjes algebras \and $QSL_p$-Spaces \and Universal representation.}

\textbf{MSC2010:} Primary 46L07; Secondary 43A30, 43A15, 47L10.
\end{abstract}

\section{Introduction}
Operator spaces were introduced in the mid 1970’s after the pioneering works
of Stinespring and Arveson, and were presented to the mathematical community
by E.G. Effros in 1986. Operator space structure on the Eymard's  Fourier-Stieltjes algebra, $B(G)$, was firstly investigated in \cite{EFFROSRUAN1991}, and fully described in \cite{EFFROSRUAN2000}. Indeed, in \cite{EFFROSRUAN2000}, operator space structure on an abstract operator space was thoroughly determined. The natural operator space structure on  Fourier-Stieltjes algebra comes from its predual $C^*$-algebra, and due to the fact that  the Fourier algebra, $A(G)$, is the predual of the von Neumann algebra $VN(G)$, generated by the left regular representation $(\lambda_{2,G},(L_2(G))$ of the locally compact group $G$, the natural operator space structure is induced on $A(G)$. Many authors benefited from this approach to investigate various problems on the Fourier and Fourier-Stieltjes algebras. For instance, in \cite{ILIESPRONK2005}, completely bounded homomorphisms on Fourier into Fourier-Stieltjes algebras are studied through the aforementioned operator space structure. 

In the next stage, by the advent of Fig\`a-Talamanca-Herz algebras $A_p(G)$, $(1<p<\infty)$, firstly by Fig\`a-Talamanca for abelian locally compact groups \cite{FIGATALAMANCA1965}, and then in general case by Herz \cite{HERZ1971}, the notion of $p$-operator space has been developed in \cite{DAWS2010}, based on the ideas of studies done by Pisier \cite{PISIER2001} and Le Merdy \cite{LEMERDY1996}. Introduced approach of $p$-operator space has been extensively utilized to turn $A_p(G)$ into a $p$-operator space, and many other properties have been studied. As a fruit of this structure on $A_p(G)$, which is a dual $p$-operator structure, Ilie has studied $p$-completely contractive homomorphisms on such algebras \cite{ILIE2013}.

 There has been a variety of approaches of defining the $p$-analog of the Fourier-Stieltjes algebras, and the most appropriate one is brought by Runde in \cite{RUNDE2005} which is denoted by $B_p(G)$. 

In this case, there can be found a need for a suitable $p$-operator space structure to be imposed. Herewith, we introduce a $p$-operator space structure on the $p$-analog of the Fourier-Stieltjes algebras, by considering it as the dual space of the $p$-operator space $UPF_p(G)$, the algebra of universal $p$-pseudofunctions. For this purpose, we divide this notes into the following sections: In Section \ref{SECTIONPRE}, as this study is based on the theory of $QSL_p$-spaces, and their $l_p$-direct sum, we briefly give essential notions and definitions about representations and their matrix form. Then, in Section \ref{SECTIONMAT}, some interesting properties of representations will be stated, and as a consequence the construction of $p$-operator space of the algebra of $p$-pseudofunctions will be demonstrated fussily. At the end, in Section \ref{SECTIONREP}, the obtained results in the latest section will be implemented to describe some features of $p$-predual space structure on $B_p(G)$.

\section{Preliminaries}\label{SECTIONPRE}
The case of group representations on Hilbert spaces  has been studied heavily. For a classic reference of primary materials on this topic, one may be referred to \cite{FOLLAND1995}, and for the operator space notions correspond to representations on Hilbert spaces one can see \cite{EFFROSRUAN2000}.  As a consequence of such attitudes, the algebra of coefficient functions of representations on Hilbert spaces, $B(G)$, was introduced in general sense by Eymard \cite{EYMARD1964}. 
The next brick in the constructing such spaces, let it be called Fourier type spaces, is the $p$-analog of the Fourier-Stieltjes algebras, and considerable number of authors have been involved to do so.
 For instance, similar constructions on semigroups has been sorted by
Amini and Medghalchi in \cite{AMINI}, and the second author in \cite{SHAMS}. Some more generalizations of
these algebras to the other settings, for examples quantum groups \cite{C-E-SCH}, and
groupoids \cite{PAT} have been studied.

  \cite{HERZ1974}, \cite{LOSERT1984}, and \cite{RUNDE2005}, took their shots at the theory of the $p$-analog of the Fourier-Stieltjes algebras, but the most suitable definition is Runde's idea, because of the fact that the introduced $B_p(G)$ coincides with Eymard's $B(G)$, whenever $p=2$. So, we follow the  definition by Runde \cite{RUNDE2005}, and to describe this $p$-analog we need the following definitions. First, we touch on the definition of specific class of Banach spaces.
   
In this section, we briefly explain necessary notions about our approach of $p$-operator space structure of the $p$-analog of the Fourier-Stieltjes algebras. We initiate this section with the definition of representations of the locally compact group $G$.
\begin{definition}
\begin{enumerate}
\item
A representation of a locally compact group $G$ is a pair of homomorphism  $\pi$ and a Banach space $E$ for which it corresponds each element of $G$ to an invertible isometric operator on $E$. Precisely, $\pi :G\rightarrow\mathcal{B}(E)$, so that $\pi(x)$ is an invertible operator with the inverse $\pi(x^{-1})$ that is isometric map.
\item
A \text{lift} of a representation of the locally compact group $G$ to the group algebra $L_1(G)$ is a contractive homomorphism $\pi :L_1(G)\rightarrow \mathcal{B}(E)$ defined through
\begin{align*}
\langle\pi(f)\xi,\eta\rangle=\int_{G}f(x)\langle\pi(x)\xi,\eta\rangle dx,\quad\xi\in E,\; \eta\in E^*.
\end{align*}
And this operator is continuous with respect to original topology on $G$ and the strong operator topology on $\mathcal{B}(E)$.
\item
A representation $(\pi,E)$ is called cyclic with the cyclic vector $\xi$, if the norm closure of the space $\pi(L_1(G))\xi$ is dense in $E$. In this case, we may denote $(\pi,E)$ by $(\pi_\xi,E_\xi)$, and call $E_\xi$ a cyclic space. 
\end{enumerate}
\end{definition}
One of the most significant associated spaces (algebras) with the representations of a locally compact group $G$, is the space of coefficient functions of them. So, we state the upcoming definition.
\begin{definition}
For a representation $(\pi, E)$ of a locally compact group $G$, we call $u$ a coefficient function of $(\pi, E)$, if there exist elements $\xi\in E$, and $\eta\in E^*$ so that
\begin{align*}
u(x)=\langle\pi(x)\xi,\eta\rangle,\quad x\in G.
\end{align*}
We may use the notation $u=u_{\xi,\eta}$ to demonstrate $u$ by emphasizing on vectors $\xi$ and $\eta$, or we may benefit from expression $u=u_{\pi,E}$ when our point is about representation $\pi$ and its associated space $E$.
\end{definition}

\begin{definition}
\begin{enumerate}
\item
We say that a Banach space $E$ is an $L_p$-space, if it is of the form of $L_p(X,\mu)$ for a measure space $(X,\mu)$.
\item
A Banach space $E$ is called a $QSL_p$-space, if it can be identified with a quotient of a subspace of an $L_p$-space.
\end{enumerate}
\end{definition}
Next definition is going to provide a relation between representations.

\begin{definition}
Let  $(\pi,E)$ and $(\rho,F)$ be two representations of a locally compact group $G$. Then
\begin{enumerate}
\item
the representations $(\pi,E)$ and $(\rho,F)$ are equivalent, if there exists an invertible isometry $T:E\rightarrow F$ such that $\rho(f)\circ T=T\circ\pi(f)$, for every $f\in L_1(G)$. In this case we write $(\pi,F)\sim(\rho,F)$.
\item
the representation $(\rho,F)$ is called a subrepresentation of $(\pi,E)$, if $F\subseteq E $, a closed subspace, and $\pi(f)=\pi(f)|_F$, for every $f\in L_1(G)$.
\item
the representation $(\rho, F)$ is said to be contained in $(\pi,E)$, and write $(\rho, F)\subset (\pi,E)$, if $(\rho, F)$ is equivalent to a subrepresentation of $(\pi,E)$.

\end{enumerate}
\end{definition}
We need to set some symbols as below.
\paragraph{Notation.}
The collection of all (classes of) representations of a locally compact group $G$ on a $QSL_p$-space is denoted by $\text{Rep}_p(G)$. Moreover, the set of all cyclic representations is denoted by $\text{Cyc}_p(G)$, as well as for a representation $(\pi,E)\in \text{Rep}_p(G)$, the set of all its cyclic subrepresentations is denoted by $\text{Cyc}_{p,\pi}(G)$.

\begin{definition}
A representation $(\pi,E)\in\text{Rep}_p(G)$ is called a $p$-universal representation, if it contains all cyclic representations in $\text{Cyc}_p(G)$.
\end{definition}
Now it is the time for introducing a type of algebras that plays a pivotal role in this paper.
\begin{definition}
Let $(\pi,E)\in\text{Rep}_p(G)$.
\begin{enumerate}
\item 
Then, if we define
\begin{align*}
\|f\|_\pi=\|\pi(f)\|,
\end{align*}
then $\|\cdot\|_\pi$ is an algebra semi-norm on $L_1(G)$.
\item
The algebra of $p$-pseudofunctions associated with $(\pi,E)$ is defined to be the operator norm closure of $\pi(L_1(G))$ in $\mathcal{B}(E)$, and is denoted by $PF_{p,\pi}(G)$. Moreover, if $(\pi,E)$ is a $p$-universal representation, then it is called the algebra of universal $p$-pseudofunctions, and we write $UPF_p(G)$ instead.
\end{enumerate}
\end{definition}
\begin{remark}\label{RUNDEREMARK}
\begin{enumerate}
\item\label{RUNDEREMARK1}
Let $(\pi,E)\in\text{Rep}_p(G)$, and put
\begin{align*}
N_\pi=\{f\in L_1(G)\ :\ \|\pi(f)\|=0\}.
\end{align*}
Then, $N_\pi$ is a closed ideal of $L_1(G)$, and by taking quotient we reach to a norm on $L_1(G)/N_\pi$. The completion of $L_1(G)/N_\pi$ with respect to this norm is $PF_{p,\pi}(G)$, described as above
\item\label{RUNDEREMARK3}
For two representations $(\pi,E)$ and $(\rho ,F)$, if $(\pi,E)$ contains all cyclic subrepresentations of $(\rho ,F)$, then for $f\in L_1(G)$, we have
\begin{align*}
\| f\|_\rho\leq\| f\|_\pi, 
\end{align*}
and this is due to the following equivalent approach of computing the norm $\|\cdot\|_\pi$,
\begin{align*}
\|\pi(f)\|=\|f\|_\pi=\sup\{\|\pi_c(f)\|=\|f\|_{\pi_c}\ :\ (\pi_c,E_c)\in\text{Cyc}_{p,\pi}(G)\}.
\end{align*}
\item\label{RUNDEREMARK4}
For $f\in L_1(G)$, we can simply consider the above norm as following:
\begin{align*}
\|\pi(f)\|=\|f\|_\pi=\sup\{\|\rho(f)\|=\|f\|_\rho\ :\ \rho\subseteq\pi\}.
\end{align*}
\item\label{RUNDEREMARK5}
Based on parts \eqref{RUNDEREMARK3} and \eqref{RUNDEREMARK4}, in the case of $UPF_p(G)$ the norm is independent of choosing specific $p$-universal representation.

\end{enumerate}
\end{remark}
In the following we define $QSL_p$-operator algebras, similar to \cite{GARDELLA2019}.

\begin{definition} Let $\mathcal{A}$ be a Banach algebra.
\begin{enumerate}
\item
A representation of $\mathcal{A}$ (on a $QSL_p$-space $E$) is a contractive homomorphism $\Pi : \mathcal{A} \rightarrow \mathcal{B}(E)$.
\item
We say that the representation $\Pi : \mathcal{A} \rightarrow \mathcal{B}(E)$ is non-degenerate (essential) if $\text{span}\{\Pi(a)\xi \ : a\in \mathcal{A},\;  \xi\in E\}$ is dense in $E$.
\item
An essential representation $\Pi : \mathcal{A} \rightarrow \mathcal{B}(E)$ of $\mathcal{A}$ (on a $QSL_p$-space $E$) is called faithful if it is injective.
\item
A  Banach algebra $\mathcal{A}$ is said to be $QSL_p$-operator algebra if there exists a $QSL_p$-space $E$ and an isometric homomorphism $\Pi : \mathcal{A} \rightarrow \mathcal{B}(E)$.
 In this case one may equivalently expect that there exists a faithful isometric representation of $\mathcal{A}$ on a $QSL_p$-space $E$.
\end{enumerate}
\end{definition}
The paper deals with $l_p$-direct sum of representations and their $(n)$-folds. So, we make a quick review on such topics. For extra information one can see \cite{DALESPOLYAKOV2012}.
\begin{remark}
For two Banach spaces $E_1$ and $E_2$, and $p\in (1,\infty )$, their $l_p$-direct sum  is the algebraic space $E_1\oplus E_2$ equipped with the following norm
\begin{align*}
\|x_1\oplus x_2\|^p=\|x_1\|^p+\|x_2\|^p, \quad x_1\oplus x_2\in E_1\oplus E_2.
\end{align*}
The completion with respect to this norm is denoted by $E_1\tilde{\oplus}_p E_2$. Additionally, we have
\begin{align*}
(E_1\tilde{\oplus}_p E_2)^*=E_1^*\tilde{\oplus}_{p'} E_2^*,\quad \frac{1}{p'}+\frac{1}{p}=1.
\end{align*}
So, for an arbitrary collection of Banach spaces $(E_i)_{i\in\mathbb{A}}$, we have
\begin{align*}
l_p-\oplus_{i\in\mathbb{A}}E_i=\{\oplus_{i\in\mathbb{A}}x_i \ :\ x_i\in E_i,\; \text{for}\; i\in \Lambda,\; \sum_{i\in\mathbb{A}}\|x_i\|^p<\infty\}.
\end{align*}
Since $\sum_{i\in\mathbb{A}}\|x_i\|^p<\infty$, then at most countably many of $x_i$ is non-zero, so every elements $x\in \oplus_i E_i$ has a countable representation $\oplus_{k=1}^\infty x_k$. Moreover, this is a well-known fact that $l_p$-direct sum of every collection of $QSL_p$-spaces is again a $QSL_p$-space.
\end{remark}
\paragraph{Notation.}
For ease of notation we use symbol $\oplus$ instead of $l_p-\oplus$, and omit the index $p$ in $\tilde{\oplus}_p$, and denote it simply via $\oplus$, if there is no ambiguity.

Now it is the time for the $p$-analog of the Fourier-Stieltjes algebras to be defined.
\begin{definition}
For $p\in(1,\infty)$, the space of all coefficient functions of representations in $\text{Rep}_p(G)$ is called $p$-analog of the Fourier-Stieltjes algebras and denoted by $B_p(G)$.
\end{definition}

\begin{remark}
\begin{enumerate}
\item
In the definition (and consequently on the related facts) of $B_p(G)$, we have swapped indexes $p$ and its conjugate number $p'$ in comparison to what expressed in \cite{RUNDE2005}.
\item
From \cite{RUNDE2005}, the norm of an element $u\in B_p(G)$ defined to be as following
\begin{align*}
\|u\|_{B_p(G)}=\inf\{\sum_k\|\xi_k\|\|\eta_k\| \ :\ u=u_{\oplus_k\xi_k,\oplus_k\eta_k}\},
\end{align*}
where the infimum is taken all over possible representations of $u$ as a coefficient function of $l_p$-direct sum of cyclic representations $(\pi_k,E_k)_{k=1}^\infty$ of $G$, i.e. 
\begin{align*}
u(x)=\sum_{k}\langle\pi_k(x)\xi_k,\eta_k\rangle=\langle\big(\oplus_{k}\pi_k)(x)\big(\oplus_{k}\xi_k\big),\oplus_{k}\eta_k\rangle,\ x\in G.
\end{align*}
Simultaneously, one may be inclined to compute the norm of an element $u\in B_p(G)$, as following
\begin{align*}
\|u\|_{B_p(G)}=\inf\{\|\xi\|\|\eta\| \ :\ u=u_{\xi,\eta}\}.
\end{align*}
Here infimum is computed throughout expressions of $u$ as a coefficient function of a $p$-universal representation, in notation we have $u(x)=\langle\pi(x)\xi,\eta\rangle$, for $x\in G$. Equipped with this norm of functions and point multiplication, the space $B_p(G)$ is a commutative unital Banach algebra. For more details, a curious reader would be referred to \cite{RUNDE2005}.
\end{enumerate}
\end{remark}
Two of the most crucial facts about $B_p(G)$ are Lemma 6.5 and Theorem 6.6. in \cite{RUNDE2005}, which are gathered in the following theorem.
\begin{theorem}\label{PROPRUNDEDUALITY}
Let $G$ be a locally compact group and $p\in (1,\infty)$, and let $(\pi,E)\in\text{Rep}_p(G)$.
\begin{enumerate}
\item\cite[Lemma 6.5]{RUNDE2005}\label{PROPRUNDEDUALITY1}
Then for each $\phi\in PF_{p,\pi}(G)^*$ there exists a unique $u\in B_p(G)$ with $\|u\|_{B_p(G)}\leq\|\phi\|_{op}$ such that
\begin{align*}
\langle \pi(f),\phi\rangle=\int_{G}f(x)u(x)dx,\quad f\in L_1(G).
\end{align*}
Moreover, if $(\pi,E)$ is a $p$-universal representation we have $\|u\|_{B_p(G)}=\|\phi\|_{op}$.
\item\cite[Theorem 6.6]{RUNDE2005}\label{PROPRUNDEDUALITY2}
\begin{enumerate}
\item
The dual space $PF_{p,\pi}(G)^*$ embeds into the algebra $B_p(G)$  contractively.
\item
The embedding $UPF_p(G)^*$ into $B_p(G)$ is an isometric isomorphism.
\end{enumerate}
\end{enumerate}
\end{theorem}

\begin{remark}
From Theorem \ref{PROPRUNDEDUALITY}, norm of elements $\pi(f)\in UPF_p(G)$ and $u\in B_p(G)$ can be computed through the following relations:
\begin{align*}
&\|\pi(f)\|_{op}=\sup\{|\langle\pi(f),v\rangle|\ :\ v\in B_p(G),\; \|v\|_{B_p(G)}\leq 1\},\\
&\|u\|_{B_p(G)}=\sup\{|\langle\pi(g),u\rangle|\ :\ \pi(g)\in UPF_p(G),\; \|\pi(g)\|_{op}\leq 1\}.
\end{align*}

\end{remark}
Now it is the turn for the definition of the $p$-operator space structure to be revealed. Our references on this topic are \cite{DAWS2010}, and \cite{LEMERDY1996}.

 For $p\in (1,\infty)$, a measure $\mu$ and a Banach space $E$, by identifying the algebraic tensor product $L_p(\mu)\otimes E$ as a subspace of $L_p(\mu,E)$ in a natural way, one may insert a norm on this algebraic tensor product, and we denote the completion of $L_p(\mu)\otimes E$ with this norm by $L_p(\mu)\tilde{\otimes}_p E$, for which we have $L_p(\mu)\tilde{\otimes}_p E=L_p(\mu, E)$ isometrically. In the following, the space $l^n_p\tilde{\otimes}_p E$ is the same space where $l_p^n$ is $\mathbb{C}^n$ equipped with $l_p$-norm. In this case $l^n_p\tilde{\otimes}_p E$ can be considered as $E^{(n)}$, the $(n)$-fold of $E$.
\begin{definition}\label{DEFpOPERATORSPDAWS}
\begin{enumerate}
\item
A concrete $p$-operator space is a closed subspace of $\mathcal{B}(E)$, for some $QSL_p$-space $E$.

 In this case for each $n\in\mathbb{N}$ one can define a norm $\|\cdot\|_n$ on $\mathbb{M}_n(X)=\mathbb{M}_n\otimes X$, by identifying $ \mathbb{M}_n(X) $ with a subspace of $ \mathcal{B}(l_p^n\otimes_p E $), where $\mathbb{M}_n$ is the space $\mathcal{B}(l_p^n)$. So, we have the family of norms $\Big(\|\cdot\|_n\Big)_{n\in\mathbb{N}}$ satisfying:

\begin{enumerate}
\item[$\mathcal{D}_\infty:$] For $u\in \mathbb{M}_n(X)$ and $v\in \mathbb{M}_m(X)$, we have that $ \|u\oplus v\|_{n+m}=\max\{\|u\|_n,\|v\|_m\} $. Here $u\oplus v\in \mathbb{M}_{n+m}(X)$, has block representation
$\begin{pmatrix}
u & 0 \\
0 & v
\end{pmatrix}. $

\item[$\mathcal{M}_p:$] For every $u\in \mathbb{M}_m(X)$ and $\alpha\in \mathbb{M}_{n,m}$, $\beta\in \mathbb{M}_{m,n}$, we have that 
$$ \|\alpha u\beta\|_n\leq\|\alpha\|_{\mathcal{B}(l^m_p,l^n_p)}\|u\|_m\|\beta\|_{\mathcal{B}(l^n_p,l^m_p)}. $$
\end{enumerate}
\item
An abstract $p$-operator space is a Banach space $X$ equipped with the family of norms $\big(\|\cdot\|_n\big)$ defined by $\mathbb{M}_n(X)$ which satisfy two axioms above.

\end{enumerate}
\end{definition}
One of the most famous application of such $p$-operator space is done by Daws for the Fig\`a-Talamanca-Herz algebras, $A_p(G)$. In \cite{DAWS2010}, there have been introduced two distinguished $p$-operator space structures on $A_p(G)$, namely, quotient and dual $p$-operator space structures. It is shown that when $G$ is amenable locally compact group, two structures coincide. There also can be found other results that seem to be precious. In the following we state some of them. These results are highly valued in our approach towards $p$-operator space structure on $B_p(G)$. Before that we give a well-known definition.
\begin{definition}
Let $X$ and $Y$ be two $p$-operator spaces, and $\Phi :X\rightarrow Y$ be a linear map. The $(n)$-fold of the map $\Phi$ can be define naturally through:
\begin{align*}
\Phi^{(n)}: \mathbb{M}_n(X)\rightarrow\mathbb{M}_n(Y),\quad \Phi^{(n)}([x_{ij}])=[\Phi(x_{ij})].
\end{align*}
and its $p$-complete norm is
\begin{align*}
\|\Phi\|_{p-cb}=\sup_{n\in\mathbb{N}}\|\Phi^{(n)}\|.
\end{align*}
Moreover, we say that $\Phi$ is $p$-completely bounded, when $\|\Phi\|_{p-cb}<\infty$, and $p$-completely contractive, if $\|\Phi\|_{p-cb}\leq 1$. Finally, it is called a $p$-completely isometric map, if for each $n\in\mathbb{N}$, the map $\Phi^{(n)}$ is an isometric map.
\end{definition}
Now we state the following obtained facts from \cite{DAWS2010}, that help us come up with an appropriate idea about $p$-operator space structure on $B_p(G)$.
\begin{theorem}\label{DAWSTHEOREMS}
Let $X$ be a $p$-operator space.
\begin{enumerate}
\item\cite[Lemma 4.2]{DAWS2010}\label{DAWSTHEOREMS1}
Let $\mu\in X^*$. Then $\mu$ is a $p$-completely isometric map from $X$ to $\mathbb{C}$, and $\|\mu\|_{p-cb}=\|\mu\|$.
\item\cite[Theorem 4.3]{DAWS2010}\label{DAWSTHEOREMS2}
There exists a $p$-completely isometry $\Phi :X^*\rightarrow \mathcal{B}(l_p(I))$ for some index set $I$.
\end{enumerate}
\end{theorem}
\begin{remark}\label{NORMCOMPUTATIONDAWS}
It is worthwhile to take notice of the fact that for a matrix $T=[T_{st}]\in\mathbb{M}_n(X)\subset \mathbb{M}_n(\mathcal{B}(E))$ norm is computed via the following formula
\begin{align*}
\|[T_{st}]\|=\sup\{\sum_{s=1}^n\|\sum_{t=1}^nT_{st}(\xi_t)\|\ :\ (\xi_t)_{t=1}^n\subset E,\; \sum_{t=1}^n\|\xi_t\|^p\leq 1\},
\end{align*}
while for a matrix $\mu=[\mu_{ij}]\in \mathbb{M}_m(X^*)$ we have the following relation to compute the norm:
\begin{align*}
\|[\mu_{ij}]\|=\sup\{\|\langle\langle [T_{st}],[\mu_{ij}]\rangle\rangle\|\ : \ n\in\mathbb{N},\; [T_{st}]\in \mathbb{M}_n(X),\;  \|[T_{st}]\|\leq 1 \},
\end{align*}
where 
\begin{align*}
\langle\langle [T_{st}],[\mu_{ij}]\rangle\rangle=\big(\langle T_{st},\mu_{ij}\rangle\big)\in\mathbb{M}_m\otimes\mathbb{M}_n=\mathbb{M}_{mn},
\end{align*}
and the scalar matrix norm $\|\langle\langle [T_{st}],[\mu_{ij}]\rangle\rangle\|$ can be computed through taking supremum on values
\begin{align*}
|\sum_{s=1}^n\sum_{i=1}^m\sum_{t=1}^n\sum_{j=1}^m\beta_{si}\langle\pi(f_{st}),u_{ij}\rangle\alpha_{tj}|,
\end{align*}
where $\alpha=[\alpha_{tj}]\in\mathbb{M}_{1,nm}$, and $\beta=[\beta_{si}]\in\mathbb{M}_{1,nm}$ with $\sum_{tj}|\alpha_{tj}|^p\leq 1$, and $\sum_{si}|\beta_{si}|^{p'}\leq 1$. Indeed, here the identification $\mathbb{M}_{1,nm}=\mathcal{B}(l_p^{nm},l_p^1)$ is applied.
\end{remark}

\section{Matrix Norm on the Algebra of $p$-Pseudofunctions}\label{SECTIONMAT}

In this section we try to define a family of norms $(\|\cdot\|_n)_{n\in\mathbb{N}}$ on the algebra of $p$-pseudofunctions, $PF_{p,\pi}(G)$, for $(\pi,E)\in\text{Rep}_p(G)$. First, we determine our setting by clarifying about elements and behavior of the structure to guarantee that the sequel $p$-operator space structure is well-defined. In the initial step, we state straightforward proposition,  that will be used implicitly in dealing with $l_p$-direct sum of representations. Next, we bring some definitions to regulate the context, and these definitions seem to be somewhat new.

\begin{proposition}\label{PROPDIRECTBLMAP}
Let $(T_i)_{i\in\Lambda}$ be a collection of bounded linear operators   $T_i:E_i\rightarrow E_i$,  $i\in\Lambda$. Then $\oplus_{i\in\Lambda}T_i:E_i\rightarrow E_i$ for $p\in (1,\infty )$ is well-defined, and we have
\begin{align*}
\|T\|^p=\sup\{\sum_{k=1}^\infty\| T_k(x_k)\|^p\ :\ (T_k)_{k=1}^\infty\subseteq (T_i)_{i\in\Lambda},\;\oplus_{k} x_k\in\oplus_{k} E_k,\; \sum_{k}\|x_k\|^p\leq 1\}.
\end{align*}
\end{proposition}

To ease of reading, we  clarify   some notations.
\paragraph{Notation.}\label{NTTN}
\begin{enumerate}
\item\label{NTTN1}
Let $(\pi,E)\in\text{Rep}_p(G)$. For $n\in\mathbb{N}$, by $\big(\pi^{(n)},E^{(n)}\big)$ we mean:
\begin{align*}
&E^{(n)}=l_p^{(n)}\tilde{\otimes}_p E;\\
&\pi^{(n)}:\mathbb{M}_n(L_1(G))\rightarrow \mathcal{B}(E^{(n)}),\quad \pi^{(n)}([f_{ij}]):=[\pi(f_{ij})]:E^{(n)}\rightarrow E^{(n)};\\
\end{align*}
where for $\; (x_j)_{j=1}^n\in E^{(n)}$, and $ [f_{ij}]\in \mathbb{M}_n(L_1(G))$, we have
\begin{align*}
\pi^{(n)}([f_{ij}])(x_j)_{j=1}^n:=[\pi(f_{ij})](x_j)_{j=1}^n=\bigg(\sum_{j=1}^n\pi(f_{ij})(x_j)\bigg)_{i=1}^n\in E^{(n)}.
\end{align*}
\item\label{NTTN3}
For two Banach spaces $E_1$ and $E_2$, we have
\begin{align*}
\big(E_1\oplus E_2\big)^{(n)}=E_1^{(n)}\oplus E_2^{(n)}.
\end{align*}
and consequently, for $T=T_1\oplus T_2: E_1\oplus E_2\rightarrow E_1\oplus E_2$ we have
\begin{align*}
T^{(n)}=\big(T_1\oplus T_2\big)^{(n)}=T_1^{(n)}\oplus T_2^{(n)}.
\end{align*}

\end{enumerate}

\begin{definition}
\begin{enumerate}
\item
For a representation $(\pi,E)\in\text{Rep}_p(G)$, we call the pair $(\pi^{(n)},E^{(n)})$ a matrix representation
of $L_1(G)$ on the $(n)$-fold of the $QSL_p$-space $E$.
\item
With the notations from above, we call $(\Theta,K)$ a cyclic matrix representation, whenever $K\subset E^{(n)}$ is a closed subspace, and $\Theta=\pi^{(n)}|_K$ is a matrix representation from $\mathbb{M}_n(L_1(G))$ into $\mathcal{B}(K)$, and there exists $\mathbf{x}\in K$ for which we have $\overline{\Theta(\mathbb{M}_n(L_1(G)))\mathbf{x}}=K$.
\end{enumerate}
\end{definition}

Next lemma illustrates the relations between $(n)$-folds of cyclic representations and cyclic matrix representations.

\begin{lemma}\label{APPLICLEMMA}
\begin{enumerate}
\item\label{APPLICLEMMA1}
Let $(\pi,E)\in\text{Rep}_p(G)$ be a representation constructed of the $l_p$-direct sum of two representations $(\pi_1,E_1)$ and  $(\pi_2,E_2)$. Then cyclic subrepresentations of $(\pi,E)$ are the $l_p$-direct sum of cyclic subrepresentations of $(\pi_1,E_1)$ and $(\pi_2,E_2)$, and vice versa.
\item\label{APPLICLEMMA2}
Let $(\pi,E)\in\text{Rep}_p(G)$. If $(\rho, F)$ is cyclic subrepresentation of $(\pi,E)$, then $(\rho^{(n)},F^{(n)})$ is a cyclic subrepresentation of $(\pi^{(n)},E^{(n)})$, that is indeed a cyclic matrix representation. Moreover, if $(\Theta ,K)$ is a cyclic matrix representation contained in $(\pi^{(n)},E^{(n)})$, then there exist cyclic subrepresentations $(\rho_k,F_k)_{k=1}^n$ that satisfy
\begin{align*}
((\oplus_{k}\rho_k)^{(n)},(\oplus_{k} F_k)^{(n)})= (\Theta ,K).
\end{align*}
\end{enumerate}
\begin{proof}
\begin{enumerate}
\item
For $k=1,2$, let $(\rho_k,F_k)\subseteq (\pi_k,E_k)$ be cyclic subrepresentations with cyclic vectors $x_k$. Then obviously $(\rho_1\oplus\rho_2,F_1\oplus F_2)$ is a cyclic subrepresentation of $(\pi,E)$ with the cyclic vector $x_1\oplus x_2$. Precisely, let $y\in F_1\oplus F_2$ be of the form $y_1\oplus y_2$. Then for a given $\varepsilon>0$ there exist $f_k\in L_1(G)$, for $k=1,2$ such that
\begin{align*}
\|\pi_1(f_1)x_1-y_1\|\leq\varepsilon,\quad \|\pi_2(f_2)x_2-y_2\|\leq \varepsilon,
\end{align*}
then we have
\begin{align*}
\|\rho_1(f_1)x_1\oplus\rho_2(f_2)x_2-(y_1\oplus y_2)\|^p=\|\pi_1(f_1)x_1-y_1\|^p+\|\pi_1(f_1)x_2-y_2\|^p\leq 2\varepsilon^p.
\end{align*}
On the other hand, we have $\rho_1(f_1)x_1\oplus\rho_2(f_2)x_2\in \rho_1\oplus\rho_2(L_1(G))(x_1\oplus x_2)$. So, the space $\rho_1\oplus\rho_2(L_1(G))(x_1\oplus x_2)$ is dense in $F_1\oplus F_2$.\\
Conversely, if $(\rho,F)\subset (\pi,E)$ is a cyclic representation with cyclic vector $x=x_1\oplus x_2\in E=E_1\oplus E_2$. Then
\begin{align*}
F=\overline{\rho(L_1(G))x}=\overline{\rho(L_1(G))x_1\oplus\rho(L_1(G))x_2},
\end{align*}
together with the fact that $\rho=\pi|_F$. Put 
\begin{align*}
F_k=\overline{\pi_k(L_1(G))x_k},\quad\text{and}\quad\rho_k=\pi_k|_{F_k},\quad \text{for}\; k=1,2.
\end{align*}
Then, for each $k=1,2$ the representation $(\rho_k,F_k)$ is a cyclic subrepresentation of $(\pi_k,E_k)$, and from previous part, $(\rho_1\oplus\rho_2,F_1\oplus F_2)$ is cyclic representation of $(\pi,E)$ so that
\begin{align*}
\rho_1\oplus\rho_2=\pi{|_{F_1\oplus F_2}}.
\end{align*}
Additionally,
\begin{align*}
F=\overline{\rho(L_1(G))(x_1\oplus x_2)}&=\overline{\pi(L_1(G))(x_1\oplus x_2)}\\
&=\overline{\pi(L_1(G))x_1}\oplus \overline{\pi(L_1(G)) x_2}\\
&=\overline{\pi_1(L_1(G))x_1}\oplus \overline{\pi_2(L_1(G)) x_2}\\
&=\overline{\rho_1(L_1(G))x_1}\oplus\overline{ \rho_2(L_1(G)) x_2}\\
&=F_1\oplus F_2,
\end{align*}
and consequently, $\rho=\rho_1\oplus\rho_2$.
\item
Let $(\rho,F)\subseteq (\pi,E)$ be cyclic subrepresentation with cyclic vector $x\in F$. Let $\mathbf{y}\in F^{(n)}$, and $\mathbf{y}=(y_1,\ldots,y_n)$. Then, for a given $\varepsilon>0$ there exist $f_i\in L_1(G)$ for $i=1,\ldots, n$, such that
\begin{align*}
\|\rho(f_i)x-y_i\|\leq \varepsilon,\quad i=1,\ldots,n.
\end{align*}
Then
\begin{align*}
\|\pi^{(n)}([F_{ij}])\mathbf{x}-\mathbf{y}\|^p=\sum_{i=1}^n\|\rho(f_i)x-y_i\|^p\leq n\varepsilon^p,
\end{align*}
 where
 \begin{align*}
 \mathbf{x}=(x,0,\ldots, 0),\quad F_{ij}=\left\{
\begin{array}{ll}
F_{i1}=f_i,&\text{if}\; j=1\\
0&\text{if}\; j\neq 1
\end{array}\right..
 \end{align*}
 Consequently, $\overline{\rho^{(n)}(L_1(G))\mathbf{x}}=F^{(n)}$.\\
 Conversely, let $(\Theta, K)$ be such that $K\subset F^{(n)}$, is a closed subspace, and $\Theta=\pi^{(n)}|_K$ has the property that there exists $\mathbf{x}\in K$, so that
 \begin{align*}
 \overline{\Theta (\mathbb{M}_n(L_1(G)))\mathbf{x}}=K,\quad \mathbf{x}=(x_1,\ldots,x_n).
 \end{align*}
 Consider the following closed subspaces:
 \begin{align*}
 F_k=\overline{\pi(L_1(G))x_k},\quad k=1,\ldots , n.
 \end{align*}
Then, put
\begin{align*}
\rho_k=\pi|_{F_k},\quad k=1,\ldots, n.
\end{align*}
Now, we claim that
\begin{align*}
K=(\oplus_{k=1}^n F_k)^{(n)},
\end{align*}
and consequently,
\begin{align*}
\Theta=(\oplus_{j=k}^n\rho_k)^{(n)}.
\end{align*}
First, we need to take notice of the structure of spaces $K$ and $\oplus_{k}F_k^{(n)}$. Since $K=\overline{\Theta(\mathbb{M}_n(L_1(G)))\mathbf{x}}$, then the space of vectors of the following form is dense in $K$:
\begin{align}\label{SPACE1}
\pi^{(n)}([f_{ij}])\mathbf{x}=\left[
\begin{array}{lll}
&\sum_{j}\pi(f_{1j})x_j&\\
&\vdots&\\
&\sum_{j}\pi(f_{ij})x_j&\\
&\vdots&\\
&\sum_{j}\pi(f_{nj})x_j&
\end{array}\right].
\end{align}
On the other hand, if for each $k=1,\ldots,n$, we set 
$$\mathbf{x}_k=(0,\ldots,0,\underbrace{x_k}_{\text{k-th}},0,\ldots,0),$$
then we have $\oplus_k\mathbf{x}_k=\mathbf{x}$, and the space $\oplus_{k}F^{(n)}_k$, the associated space with cyclic representation $\oplus_{k}\rho^{(n)}_k$, is generated by cyclic vector $\oplus\mathbf{x}_k$,  so $\oplus_kF^{(n)}_k$ is the closure of the space of vectors of the form
\begin{align*}
&(\oplus_k\rho_k)^{(n)}([g_{ij}])\oplus_k\mathbf{x}_k=\oplus_k([\rho_k(g_{ij})]\mathbf{x}_k),
\end{align*}
with matrix representation as following:
\begin{align}\label{SPACE2}
&\oplus_k\left[
\begin{array}{l}
\rho_k(g_{1k}) x_k\\
\vdots\\
\rho_k(g_{ik}) x_k\\
\vdots\\
\rho_k(g_{nk}) x_k
\end{array}\right]
=\left[
\begin{array}{l}
\oplus_k\rho_k(g_{1k}) x_k\\
\vdots\\
\oplus_k\rho_j(g_{ik}) x_k\\
\vdots\\
\oplus_k\rho_k(g_{nk}) x_k
\end{array}\right]
=\left[
\begin{array}{l}
\sum_k\pi(g_{1k}) x_k\\
\vdots\\
\sum_k\pi(g_{ik}) x_k\\
\vdots\\
\sum_k\pi(g_{nk}) x_k
\end{array}\right],
\end{align}
 In above, we have benefited from the fact that for each $k=1,\ldots,n$, the representation $\rho_k$ is the restriction of $\pi$ to $F_k$. So, as it is clear from the appearance of vectors in \eqref{SPACE1} and \eqref{SPACE2}, we have $K=\oplus_k F^{(n)}_k$, and therefore,
\begin{align*}
\Theta=\pi^{(n)}|_K=\pi^{(n)}|_{\oplus_kF^{(n)}_k}=\oplus_k\rho^{(n)}_k.
\end{align*}

\end{enumerate}

\end{proof}
\end{lemma}

\begin{remark}
As a summary of the previous lemma, we have one to one correspondence between $(n)$-folded cyclic subrepresentations and cyclic matrix representations. So, our new definition of cyclic matrix representations is nothing except rephrasing of the $(n)$-folds of $l_p$-direct sum of $n$ cyclic subrepresentations. The Lemma \ref{APPLICLEMMA} has an implicitly intricate role as we shall work with the norm of matrix representations.
\end{remark}

Next lemma is necessary at the aim of clarifying the matrix spaces of the algebra of $p$-pseudofunctions.

\begin{lemma}\label{LemmaMatrixPSE}
 Let $(\pi,E)\in\text{Rep}_p(G)$. Then, we have
\begin{align*}
\overline{\pi^{(n)}(\mathbb{M}_n(L_1(G))){}}^{\|\cdot\|_{\mathcal{B}(E^{(n)})}}=\mathbb{M}_n(PF_{p,\pi}(G)).
\end{align*}
\begin{proof}

Let $T=[T_{ij}]\in \mathbb{M}_n(PF_{p,\pi}(G))$. Then for each pair $(i,j)$ there exists a net $(\pi(f_{\alpha_{ij}}))_{\alpha_{ij}\in\Lambda_{ij}}$, where $\Lambda_{ij}$ is a directed set, such that
\begin{align*}
T_{ij}=\|\cdot\|-\lim_{\alpha_{ij}\in\Lambda_{ij}}\pi(f_{\alpha_{ij}}).
\end{align*}
Now, put
\begin{align*}
\Lambda=\Pi_{i=1}^n\Pi_{j=1}^n\Lambda_{ij},
\end{align*}
then $\Lambda$ is a directed set  as well, through 
\begin{align*}
\alpha\leq \beta\Longleftrightarrow \alpha_{ij}\leq \beta_{ij},
\end{align*}
where $\alpha=[\alpha_{ij}]\in \Lambda$, and $\beta=[\beta_{ij}]\in \Lambda$. Let
\begin{align*}
[\pi^{(n)}(F_\alpha)]=[\pi (f_{\alpha_{ij}})].
\end{align*}
Consequently, we have
\begin{align*}
[T_{ij}]=\|\cdot\|_n-\lim_{\alpha}[\pi^{(n)}(F_\alpha)],
\end{align*}
and this is due to the fact that for any matrix $[x_{ij}]$, we have
\begin{align*}
\|[x_{ij}]\|^p\leq\sum_{i,j}\|x_{ij}\|^p.
\end{align*}
More precisely, for a given $\varepsilon>0$, let $\alpha^0_{ij}\in\Lambda_{ij}$ be such that
\begin{align*}
\|T_{ij}-\pi(f_{\alpha_{ij}})\|<\varepsilon,\qquad \text{for}\ \alpha_{ij}\geq \alpha^0_{ij}.
\end{align*}
Then we have
\begin{align*}
\|[T_{ij}-\pi(f_{\alpha_{ij}})]\|^p\leq\sum_{i,j}\|T_{ij}-\pi(f_{\alpha_{ij}})\|^p<n^2\varepsilon^p,
\end{align*}
for every $\Lambda\ni\alpha=[\alpha_{ij}]\geq[\alpha^0_{ij}]=\alpha^0$.
\end{proof}

\end{lemma}

\begin{remark}
 For $(\pi,E)\in\text{Rep}_p(G)$. Then the Lemma has clarified the following formula
 \begin{align*}
\overline{\pi^{(n)}(\mathbb{M}_n(L_1(G))){}}^{\|\cdot\|_{\mathcal{B}(E^{(n)})}}=\mathbb{M}_n(PF_{p,\pi}(G))=\mathbb{M}_n\big({}\overline{(\pi(L_1(G)){}}^{\|\cdot\|_{\mathcal{B}(E)}}\big),
\end{align*}
which is a passage of closure of the norm $\|\cdot\|_{\mathcal{B}(E)}$ to the matrix forms. 
\end{remark}

The sequel proposition is our first step toward defining a $p$-operator space structure on the space $B_p(G)$. First, we provide evidences indicating that the algebra $PF_{p,\pi}(G)$ is a $QSL_p$-operator algebra.

\begin{proposition}\label{PROPQSLPALGEBRA}
Let $(\pi,E)\in\text{Rep}_p(G)$. Then the Banach algebra $PF_{p,\pi}(G)$ of $p$-pseudofunctions associated with $(\pi,E)$, is a $QSL_p$-operator algebra.
\begin{proof}
To prove the claim we need to find an isometric representation of $PF_{p,\pi}(G)$. Let $g\in L_1(G)$ and $r\in\mathbb{N}$. From Remark \ref{RUNDEREMARK}, there exists a cyclic subrepresentation $(\pi_{g,r},E_{g,r})$ of $(\pi,E)$ with cyclic vector $\xi_{g,r}$ such that
\begin{align*}
\|\pi_{g,r}(g)\xi_{g,r}\|>\|\pi(g)\|-\frac{1}{r},\quad \|\xi_{g,r}\|\leq 1.
\end{align*}
Consider the following $l_p$-direct sum of aforementioned cyclic subrepresentations associated with the $l_p$-direct sum of their paired spaces.
\begin{align}\label{FIRSTDEFPI}
&\mathcal{E}=\oplus_{g\in L_1(G)}\oplus_{r=1}^\infty E_{g,r}\qquad\Pi=\oplus_{g\in L_1(G)}\oplus_{r=1}^\infty\pi_{g,r}.
\end{align}
Each element $x\in\mathcal{E}$ has countable representation, that is of the form
\begin{align*}
x=\oplus_{k=1}^\infty\oplus_{r=1}^\infty x_{g_k,r},
\end{align*}
and we have
\begin{align*}
\|x\|^p=\sum_{k=1}^\infty\sum_{r=1}^\infty\|x_{g_k,r}\|^p.
\end{align*}
Now, for an arbitrary $f\in L_1(G)$, we calculate the norm $\|\Pi(f)\|$. We have
\begin{align*}
\Pi(f)=\oplus_{g\in L_1(G)}\oplus_{r=1}^\infty\pi_{g,r}(f).
\end{align*}
Moreover, from Proposition \ref{PROPDIRECTBLMAP}, in order to compute the norm $\|\Pi(f)\|$, we have the following formula:
\begin{align}\label{NORMPI}
\|\Pi(f)\|^p=\sup\{\sum_{k,r}\|\pi_{g_k,r}(f)x_{g_k,r}\|^p\ :\  x=\oplus_{k,r}x_{g_k,r}\in\mathcal{E},\;\sum_{k,r}\|x_{g_k,r}\|^p\leq 1\}.
\end{align}
In the rest of the proof we shall show that $(\Pi, \mathcal{E})$ is an isometric representation of $PF_{p,\pi}(G)$ into a closed subalgebra of $\mathcal{B(E)}$. For this aim, first we deal with the dense subspace $\pi(L_1(G))$. Let $f\in L_1(G)$ and $\varepsilon>0$ is given. Then via \eqref{NORMPI}, there exists $x=\oplus_{k,r} x_{g_k,r}\in\mathcal{E}$ such that
\begin{align*}
\|\Pi(f)\|^p-\varepsilon <\sum_{k,r}\|\pi_{g_k,r}(f)x_{g_k,r}\|^p,\quad \text{and},\quad \sum_{k,r}\|x_{g_k,r}\|^p\leq 1.
\end{align*}
On the other hand, each $\pi_{g_k,r}$ is a restriction of $\pi$ to the subspace $E_{g_k,r}$, and it implies that
\begin{align*}
\pi_{g_k,r}(f)x_{g_k,r}=\pi(f)x_{g_k,r},\quad \text{and},\quad\|\pi_{g_k,r}(f)\|\leq \|\pi(f)\|.
\end{align*}
Utilizing this fact obtains
\begin{align*}
\|\Pi(f)\|^p-\varepsilon &<\sum_{k,r}\|\pi_{g_k,r}(f)x_{g_k,r}\|^p\\
&\leq\sum_{k,r}\|\pi_{g_k,r}(f)\|^p \|x_{g_k,r}\|^p\\
&\leq\sum_{k,r}\|\pi(f)\|^p \|x_{g_k,r}\|^p\\
&\leq \|\pi(f)\|^p\sum_{k,r}\|x_{g_k,r}\|^p\\
&\leq \|\pi(f)\|^p.
\end{align*}
Consequently, since $\varepsilon>0$ is arbitrary, we have $\|\Pi(f)\|\leq\|\pi(f)\|$.

We need to take notice of the fact that since $\pi$ is contractive homomorphism, last inequality implies that $\Pi$ is as well.\\
Additionally, for the inverse inequality, if by \eqref{FIRSTDEFPI}, for $m\in\mathbb{N}$,  we choose $x_m=\oplus_{k,r} x^{(m)}_{g_k,r}\in\mathcal{E}$ such that
\begin{align}\label{TECHEQ}
x^{(m)}_{g_k,r}=\left\{
\begin{array}{ll}
x^{(m)}_{g_1,m}=\xi_{f,m}& \text{if}\; k=1,\; r=m\\
0&\text{o.w.}
\end{array}\right.,
\end{align}
then we have $\|x_m\|\leq 1$, and
\begin{align*}
\|\Pi(f)\|\geq\|\Pi(f)x_m\|=\|\pi_{f,m}(f)\xi_{f,m}\|>\|\pi(f)\|-\frac{1}{m},
\end{align*}
which means $\|\Pi(f)\|\geq \|\pi(f)\|$.

So, $\Pi : \pi(L_1(G))\rightarrow\mathcal{B(E)}$ is an isometric representation, that can be extended to $PF_{p,\pi}(G)$ isometrically.
\end{proof}

\end{proposition}

\begin{remark}\label{REMPROPQSLPALGEBRA}
 
The idea of this proof is come from Proposition 4.6 of \cite{GARDELLA2019}. However, here we use $l_p$-direct sum of cyclic representations rather than $l_\infty$-one. Besides, here we gave a general statement about the algebras of $p$-pseudofunctions not just for the universal one.

\end{remark}

\begin{lemma}\label{COMPUTLEMMA}
Let the representation $(\Pi,\mathcal{E})$ be as mentioned above. Then we have
\begin{align*}
&(\oplus_{g\in L_1(G)}\oplus_{r=1}^\infty \pi_{g,r})^{(n)}=\Pi^{(n)}=\oplus_{g\in L_1(G)}\oplus_{r=1}^\infty \pi_{g,r}^{(n)},\\ &(\oplus_{g\in L_1(G)}\oplus_{r=1}^\infty E_{g,r})^{(n)}=\mathcal{E}^{(n)}=\oplus_{g\in L_1(G)}\oplus_{r=1}^\infty E_{g,r}^{(n)}.
\end{align*}
\begin{proof}
The claim is true due to some obvious modifications in the computations in Notation \ref{NTTN}-\eqref{NTTN3}, first for countable and then for uncountable one. We give a quick sketch of everything. First we take care of the vectors in $\mathcal{E}^{(n)}$. Let
\begin{align*}
\mathbf{x}=(x_1,\ldots,x_n)\in\mathcal{E}^{(n)}=\bigg(\oplus_{g\in L_1(G)}\oplus_{r=1}^\infty E_{g,r}\bigg)^{(n)}.
\end{align*}
For each $j=1,\ldots, n$, let $\oplus_{k,r} x^{j}_{g_k,r}\in \oplus_{k,r} E^j_{g_k,r}$ be the countable representation of $x_j$, then by adding some zero component to these representations, and rearrange them, we can assume that $\oplus_{k,r}E^j_{g_k,r}=\oplus_{k,r}E_{g_k,r}$, and therefore, the associated representation is $\oplus_{k,r}\pi_{g_k,r}$, instead of $n$ distinguished countably infinite $l_p$-direct sums $\oplus_{k,r}\pi^j_{g_k,r}$, for $j=1,\ldots,n$.
So, by such modifications in the appearance of the components $x_j$, together with the fact that we have
\begin{align}\label{NFOLDVECTOR}
\big(x_1,\ldots,x_n\big)&=\bigg(\oplus_{k,r} x^{1}_{g_k,r},\ldots,\oplus_{k,r} x^{n}_{g_k,r}\bigg)=\oplus_{k,r}\big( x^{1}_{g_k,r},\ldots, x^{n}_{g_k,r}\big),
\end{align}
then $(\oplus_{k,r}E_{g_k,r})^{(n)}=\oplus_{k,r}E_{g_k,r}^{(n)}$, so we have $\mathcal{E}^{(n)}=\oplus_{g,r}E^{(n)}_{g,r}$. Additionally, for a countable representation of $\mathbf{x}$, we have
\begin{align}\label{NFOLDMATRIX}
\Pi^{(n)}([f_{ij}])\mathbf{x}=(\oplus_{k,r}\pi_{g_k,r})^{(n)}([f_{ij}])\mathbf{x}=\oplus_{k,r}\pi_{g_k,r}^{(n)}([f_{ij}])\mathbf{x},
\end{align}
which means $(\oplus_{k,r}\pi_{g_k,r})^{(n)}=\oplus_{k,r}\pi_{g_k,r}^{(n)}$, and consequently, $\Pi^{(n)}=\oplus_{g,r}\pi_{g,r}^{(n)}$.
\end{proof}
\end{lemma}

\begin{proposition}\label{PROPPOPERATORALGRBRA}
Let $(\pi,E)\in\text{Rep}_p(G)$, and let $(\Pi,\mathcal{E})$ be the associated isometric representation of $PF_{p,\pi}(G)$. Then the matrix representation $(\Pi^{(n)},\mathcal{E}^{(n)})$ is an isometric map from $\mathbb{M}_n(PF_{p,\pi}(G))$ onto a closed subspace of  $\mathbb{M}_n(\mathcal{B}(\mathcal{E}))=\mathcal{B}(\mathcal{E}^{(n)})$.
\begin{proof}
For $n\in\mathbb{N}$, consider the following $(n)$-folded map on the subspace $\mathbb{M}_n(L_1(G))$ that is dense in $\mathbb{M}_n(PF_{p,\pi}(G))$ through Lemma \ref{LemmaMatrixPSE}.
\begin{align*}
\Pi^{(n)}: \pi^{(n)}(\mathbb{M}_n(L_1(G)))\rightarrow \mathbb{M}_n(\mathcal{B(E)})=\mathcal{B}(\mathcal{E}^{(n)}),\quad \Pi^{(n)}([\pi(f_{ij})])=[\Pi(f_{ij})].
\end{align*}
We shall show that $\|[\pi(f_{ij})]\|=\|[\Pi(f_{ij})]\|$. First we note that for a matrix $[f_{ij}]$ the norm $\|[\Pi(f_{ij})]\|$ can be calculated through the explanations in the proof of Lemma \ref{COMPUTLEMMA} and Proposition \ref{PROPDIRECTBLMAP}, by the following relation
\begin{align}\label{NFOLDNORMPI}
\|[\Pi(f_{ij})]\|^p=\sup\{\sum_{k,r}\sum_{i=1}^n\bigg\|\sum_{j=1}^n \pi_{g_k,r}(f_{ij})x^j_{g_k,r}\bigg\|^p\ :\  \sum_{j=1}^n\sum_{k,r} \| x^{j}_{g_k,r}\|^p\leq 1\}.
\end{align}
Therefore, for a given $\varepsilon>0$, we may assume that there exists a vector $\mathbf{x}\in \mathcal{E}^{(n)}$ of the form \eqref{NFOLDVECTOR} such that
\begin{align*}
\|\mathbf{x}\|^p=\sum_{j=1}^n\|x_j\|^p<\sum_{j=1}^n\sum_{k,r}\| x^{j}_{g_k,r}\|^p+\varepsilon,
\end{align*}
with $\|\mathbf{x}\|\leq 1$, and via \eqref{NFOLDNORMPI}, we have
\begin{align*}
\|[\Pi(f_{ij})]\|^p-\varepsilon &<\sum_{k,r}\sum_{i=1}^n\bigg\|\sum_{j=1}^n \pi_{g_k,r}(f_{ij})x^j_{g_k,r}\bigg\|^p .
\end{align*}
On the other hand, since for each $r$, $k$, and $f_{ij}$, the operator $\pi_{g_k,r}(f_{ij})$ is the restriction of $\pi(f_{ij})$ to the subspace ${E_{g_k,r}}$, then we have 
\begin{align*}
\pi_{g_k,r}(f_{ij})x^j_{g_k,r}=\pi(f_{ij})x^j_{g_k,r}.
\end{align*}
If we put
\begin{align*}
&y_j=\oplus_{k,r} x^j_{k,r},\quad j=1,\ldots, n,
\end{align*}
as well as
\begin{align*}
&\mathbf{y}=(y_1,\ldots, y_n),
\end{align*}
which are allowed to do so because of the Proposition \ref{PROPQSLPALGEBRA},  then we have
\begin{align*}
\|\mathbf{y}\|^p=\sum_{j=1}^n\|y_j\|^p=\sum_{k=1}^\infty\sum_{r=1}^\infty\sum_{j=1}^n\|x^j_{g_k,r}\|^p\leq 1,
\end{align*}
such that
\begin{align*}
\|[\Pi(f_{ij})]\|^p-\varepsilon &<\sum_{k=1}^\infty\sum_{r=1}^\infty\bigg(\sum_{i=1}^n\bigg\|\sum_{j=1}^n \pi^j_{g_k,r}(f_{ij})x^j_{g_k,r}\bigg\|^p\bigg)\\
&=\sum_{i=1}^n\|\sum_{j=1}^n \pi(f_{ij}) y_j\|^p,
\end{align*}
and it can be concluded that $\|[\Pi(f_{ij})]\|\leq\|[\pi(f_{ij})]\|$.

Conversely, from the definition of the norm of an element $[\pi(f_{ij})]$, for an arbitrary $m\in \mathbb{N}$, there exists a vector
\begin{align*}
\mathbf{y}^m=(y^m_j)_{j=1}^n=(y^m_1,\ldots, y^m_n)\in E^{(n)},
\end{align*}
with
\begin{align*}
\|\mathbf{y}^m\|^p=\|(y^m_j)_{j=1}^n\|^p=\|(y^m_1,\ldots, y^m_n)\|^p=\sum_{j=1}^n\|y^m_j\|^p\leq 1,
\end{align*}
such that
\begin{align*}
\|[\pi(f_{ij})]\|^p-\frac{1}{m}<\sum_{i=1}^n\|\sum_{j=1}^n\pi(f_{ij})y^m_j\|^p.
\end{align*}
Now, similar to what we have done for the case $n=1$ in \eqref{TECHEQ}, for this given and fixed $m\in\mathbb{N}$, we construct a vector $\mathbf{x^m}\in\mathcal{E}^{(n)}$ as following:
\begin{align*}
\mathbf{x}^m= \big(x^m_1,\ldots,x^m_n\big)&=\bigg(\oplus_{k,r} x^{m,1}_{g_k,r},\ldots,\oplus_{k,r} x^{m,n}_{g_k,r}\bigg)=\oplus_{k,r}\big( x^{m,1}_{g_k,r},\ldots x^{m,n}_{g_k,r}\big).
\end{align*}
 For each $j=1,\ldots, n$, put
\begin{align*}
x^{m,j}_{g_k,r}=\left\{
\begin{array}{ll}
x^{m,j}_{g_1,m}=y^m_j& \text{if}\; k=1,\; r=m\\
0&\text{o.w.}
\end{array}\right.,
\end{align*}
then we have
\begin{align*}
\|\mathbf{x}^m\|^p=\sum_{j=1}^n\sum_{k=1}^\infty\sum_{r=1}^\infty \| x^{m,j}_{g_k,r}\|^p=\sum_{j=1}^n\|y^m_j\|^p\leq 1.
\end{align*}
On top of that, we have
\begin{align}\label{THEREL}
\|[\Pi(f_{ij})]\|^p\geq\sum_{k,r}\sum_{i=1}^n\bigg\|\sum_{j=1}^n\pi^{m,j}_{g_k,r}(f_{ij})x^{m,j}_{g_k,r}\bigg\|^p=\sum_{i=1}^n\|\sum_{j=1}^n\pi(f_{ij})y^{m}_{j}\|^p>\|[\pi(f_{ij})]\|^p-\frac{1}{m},
\end{align}
which means $\|[\Pi(f_{ij})]\|\geq\|[\pi(f_{ij})]\|$. So the map $\Pi^{(n)} : \mathbb{M}_n(\pi(L_1(G)))\rightarrow\mathcal{B}(\mathcal{E}^{(n)})$ is an isometric representation, that can be extended to $\mathbb{M}_n(PF_{p,\pi}(G))$ isometrically via Lemma \ref{LemmaMatrixPSE}, and similar to the case $n=1$ in Proposition \ref{PROPQSLPALGEBRA}.
\end{proof}
\end{proposition}

\begin{theorem}\label{THEOREMRESULTSEC2}
For a representation $(\pi,E)\in\text{Rep}_p(G)$, the algebra of $p$-pseudofunctions $PF_{p,\pi}(G)$ is a $p$-operator space.
\begin{proof}
From Propositions \ref{PROPQSLPALGEBRA} and \ref{PROPPOPERATORALGRBRA}, we know that $PF_{p,\pi}(G)$ is $p$-completely isometric to a subspace of $\mathcal{B(E)}$, and due to the Definition \ref{DEFpOPERATORSPDAWS} it is an abstract $p$-operator space.
\end{proof}
\end{theorem}

\section{$p$-Operator Space Structure of $B_p(G)$}\label{SECTIONREP}

In this section, first we illustrate the relation between two algebras of $p$-pseudofunctions of two related representations in $\text{Rep}_p(G)$. Next, we prove that the $(n)$-matrix norms of $p$-universal representations are equal, and as a subsequence, the algebra of universal $p$-pseudofunctions is independent of choosing $p$-universal representation in the matrix form. Finally, it is obtained that $p$-operator space structure on $B_p(G)$ is well-defined. \begin{proposition}
If two representations $(\pi,E),(\rho, F)\in\text{Rep}_p(G)$ are such that $(\pi,E)\sim(\rho, F)$ via the map $T:E\rightarrow F$, then $(\pi^{(n)},E^{(n)})\sim (\rho^{(n)},F^{(n)})$ through the invertible isometric map $T^{(n)}:E^{(n)}\rightarrow F^{(n)}$ which is the $(n)$-fold of $T$. Consequently, two algebras $PF_{p,\pi}(G)$ and $PF_{p,\rho}(G)$ are $p$-completely isometrically isomorphism.
\begin{proof}
It is obvious fact that $T^{(n)}:E^{(n)}\rightarrow F^{(n)}$ is an invertible isometric map with the inverse $(T^{(n)})^{-1}=(T^{-1})^{(n)}$. Let us define
\begin{align*}
\Phi:PF_{p,\pi}(G)\rightarrow PF_{p,\rho}(G),\quad \Phi(\pi(f))=T\circ\pi(f)\circ T^{-1}=\rho(f),\quad f\in L_1(G).
\end{align*}
Since $T$ is an invertible isometric map then $\Phi$ is an isometric isomorphism. Now, consider the $(n)$-fold map $\Phi^{(n)}$:
\begin{align*}
\Phi^{(n)}:\mathbb{M}_n(PF_{p,\pi}(G))\rightarrow \mathbb{M}_n(PF_{p,\rho}(G)).
\end{align*}
The same is true since $T^{(n)}$ is an invertible isometric map.
\end{proof}
\end{proposition}
Next lemma is provided to assure that every cyclic representation behaves as a cyclic representation when it is contained in an arbitrary representation.
\begin{lemma}
If $(\pi_\xi,E_\xi)$ is a cyclic representation and $(\pi_\xi,E_\xi)\sim (\rho,F)$, then $(\rho, F)$ is cyclic as well.
\begin{proof}
Since $(\pi_\xi,E_\xi)\sim (\rho,F)$, then there exists an invertible isometric map $T:E_\xi\rightarrow F$ such that for every $f\in L_1(G)$, we have $T\circ\pi_\xi(f)=\rho(f)\circ T$. Put: $F\ni\phi=T\xi=T(\lim_{i\in\Lambda}\pi_\xi(e_i)\xi)$, where $(e_i)_{i\in\Lambda}\subset L_1(G)$ is bounded approximate identity. Let $\phi'\in F$ be an arbitrary vector. Since $T$ is an invertible isometry then there is a unique $\xi'\in E_\xi$ such that $T\xi'=\phi'$. Since $\xi'\in E_\xi$, then there exists a net $(\pi_\xi(f_i)\xi)_{i\in\Lambda}$ such that $\xi'=\|\cdot\|_{E_\xi}-\lim_{i\in\Lambda}\pi_\xi(f_i)\xi$. We have:
\begin{align*}
\phi'=T\xi'&=T(\lim_{i\in\Lambda}\pi_\xi(f_i)\xi)\\
&=\lim_{i\in\Lambda}T(\pi_\xi(f_i)\xi)\\
&=\lim_{i\in\Lambda}\rho(f_i)T\xi\\
&=\lim_{i\in\Lambda}\rho(f_i)\phi.
\end{align*}
and it means that each element $\phi'\in F$ can be approximated by the element $\phi=T\xi$. Consequently, $F=F_\phi$, and then $\rho=\rho_\phi$.
\end{proof}
\end{lemma}
Now, as a fruit of previous statements we have the following precious corollary.
\begin{corollary}\label{CORCYC}
If $(\pi, E)$ contains $(\rho, F)$, then $(\pi,E)$ contains all cyclic subrepresentation of $(\rho,F)$ as a cyclic subrepresentation.
\end{corollary}

\begin{proposition}\label{IMPPROP}
If $(\rho,F)$ is a subrepresentation of $(\pi,E)$, then the natural map
\begin{align*}
\Psi :PF_{p,\pi}(G)\rightarrow PF_{p,\rho}(G), \quad\pi(f)\to\rho(f),
\end{align*}
$f\in L^1(G)$, is $p$-completely contractive.
\begin{proof}
Firstly, it is clear that this map is contractive. On the other hand if we let $(\Pi,\mathcal{E})$ and $(\Theta,\mathcal{F})$ be corresponding representations to $(\pi, E)$ and $(\rho,F)$ respectively, as we had in \eqref{FIRSTDEFPI}, then we obviously have
\begin{align*}
\mathcal{F}\subseteq\mathcal{E}\quad \Theta=\Pi|_{\mathcal{F}}.
\end{align*}
Because from Corollary \ref{CORCYC} the representation $(\pi,E)$  contains all cyclic subrepresentations of $(\rho,F)$, which implies the claim.
\end{proof}
\end{proposition}

\begin{theorem}\label{UNIVERSALALGEBRA}
The algebra of universal $p$-pseudofunctions $UPF_p(G)$ is an abstract $p$-operator space and is independent of choosing specific universal representation.
\begin{proof}
Being a $p$-operator space is clear from Theorem \ref{THEOREMRESULTSEC2}, and Proposition \ref{IMPPROP} implies the dependency part.
\end{proof}
\end{theorem}
As a consequence of previous theorem, we give an immensely important theorem below.
\begin{theorem}
For $p\in (1,\infty)$, the Banach algebra $B_p(G)$ is a $p$-operator space.
\begin{proof}
By Theorem \ref{UNIVERSALALGEBRA}, the algebra $UPF_p(G)$ is a $p$-operator space. Therefore, from the combination of Theorem \ref{DAWSTHEOREMS}-\eqref{DAWSTHEOREMS2}, with Theorem \ref{PROPRUNDEDUALITY}-\eqref{PROPRUNDEDUALITY2} we have the claim.
\end{proof}
\end{theorem}
\paragraph{Notation.}\label{NTTNNOPSTNORM}
\begin{enumerate}
\item\label{NTTNNOPSTNORM1}
Let $(\pi,E)\in \text{Rep}_p(G)$ be a universal representation. Then from Theorem \ref{PROPRUNDEDUALITY}-\eqref{PROPRUNDEDUALITY1}, the duality between $UPF_p(G)$ and $B_p(G)$ is as following
\begin{align*}
\langle\pi(f),u\rangle=\int_G f(x)u(x)dx=\langle\pi(f)\xi_u,\eta_u\rangle,\quad f\in L_1(G),\; u\in B_p(G).
\end{align*}
\item\label{NTTNNOPSTNORM2}
For matrices $U=[u_{ij}]\in \mathbb{M}_m(B_p(G))$ and $F=[\pi(f_{st})]\in \mathbb{M}_n(UPF_p(G))$, from Remark \ref{NORMCOMPUTATIONDAWS}, we have
\begin{align*}
\langle\langle [\pi(f_{st})], [u_{ij}]\rangle\rangle=\langle\langle F,U\rangle\rangle:=\big(\langle\pi(f_{st}),u_{ij}\rangle\big)_{(s,i),(t,j)}\in\mathbb{M}_m\otimes\mathbb{M}_n=\mathbb{M}_{m\times n},
\end{align*}
and the norm of the element $U=[u_{ij}]\in \mathbb{M}_m(B_p(G))$ can be computed via the relation below
\begin{align*}
\|U\|&=\sup\{|\langle\langle [\pi(f_{st})],[u_{ij}]\rangle\rangle|\ :\ n\in\mathbb{N},\; [\pi(f_{st})]\in \mathbb{M}_n(UPF_p(G)),\; \|[\pi(f_{st})]\|=1\}\\
&=\sup\{|\langle\langle F,U\rangle\rangle|\ :\ n\in\mathbb{N},\; F\in \mathbb{M}_n(UPF_p(G)),\; \|F\|=1\}.
\end{align*}
Moreover, for $u\in B_p(G)$, by $\mathbf{u}_n$ we mean
\begin{align*}
\mathbf{u}_n:\mathbb{M}_n(UPF_p(G))\rightarrow \mathcal{B}(l_p^n),\quad\mathbf{u}_n([\pi(f_{st}])=\big(\langle\pi(f_{st}),u\rangle\big),
\end{align*}
and denote by $\mathbf{u}$ the function $u$ as a $p$-completely bounded operator on $UPF_{p}(G)$. Besides, from Theorem \ref{DAWSTHEOREMS}-\eqref{DAWSTHEOREMS1} we have $\|\mathbf{u}\|_{p\text{-cb}}=\|u\|_{B_p(G)}$.
\item\label{NTTNNOPSTNORM3}
The norm of an element $F=[\pi(f_{st})]\in \mathbb{M}_n({UPF}_p(G))$ can be computed via the following relation
\begin{align*}
\|F\|&=\sup\{\bigg(\sum_{s=1}^n\|\sum_{t=1}^n\langle\pi(f_{st})\xi_t\|^p\bigg)^\frac{1}{p}\ :\ (\xi_t)_{t=1}^n\subset E, \; \sum_{t=1}^n\|\xi_t\|^p\leq 1\}\\
&=\sup\{|\sum_{s=1}^n\sum_{t=1}^n\langle\pi(f_{st})\xi_t,\eta_s\rangle|\; :\; \sum_{t=1}^n\|\xi_t\|^p\leq 1,\;\sum_{s=1}^n\|\eta_s\|^{p'}\leq 1\},
\end{align*}
where $ (\xi_t)_{t=1}^n\subset E,$ and  $(\eta_s)_{s=1}^n\subset E^*$. Additionally, for $\pi(f)\in UPF_p(G)$, by $\mathbf{\pi(f)}_n$ we mean
\begin{align*}
\mathbf{\pi(f)}_n:E^{(n)}\rightarrow E^{(n)},\quad\mathbf{\pi(f)}_n((\xi_{j=1})_{j=1}^n)=\big(\pi(f)\xi_1,\ldots,\pi(f)\xi_n\big),
\end{align*}
and denote by $\mathbf{\pi(f)}$ the operator $\pi_\infty(f)$, which belongs to $UPF_{p}(G)$. In details, we have $\|\mathbf{\pi(f)}_n\|_n\leq\|\pi(f)\|$, and consequently, 
\begin{align*}
\|\mathbf{\pi(f)}\|=\|\pi_\infty(f)\|\leq\|\pi(f)\|.
\end{align*}
\end{enumerate}
     
Next proposition is the output of all materials mentioned before.
\begin{proposition}
For a locally compact group $G$, and a complex number $p\in (1,\infty)$, the identification $B_p(G)=UPF_p(G)^*$ is $p$-completely isometric isomorphism.
\end{proposition}
We believe that this result  gives an interesting point of view to study so many problems around the matrix structure of Fourier type algebras. For instance, a generalization of Ilie's researches seem to be appealing to be studied in this context, as well as a comparison of probably different $p$-operator space structures on the $p$-analog of the Fourier-Stieltjes algebras would be enticing.
\paragraph{Acknowledgements}
We would like to thank Eusebio Gardella for his great hints and providing some of his unpublished papers for us.


\begin{thebibliography}{8}



























\bibitem{AMINI}
M. Amini, A. R. Medghalchi,
\newblock \textit{Restricted algebras on inverse semigroups I, representation theory},
\newblock Math. Nachr., vol. 279, \textbf{16}, (2006), 1739-1748.

\bibitem{DALESPOLYAKOV2012}
H. G. Dales, M. E. Polyakov
\newblock   \textit{Multi-normed spaces},
\newblock   ArXiv: 1112.15148v2.












\bibitem{DAWS2010}
M. Daws,
\newblock   \textit{$p$-Operator spaces and Fig\`a-Talamanca-Herz algebras},
\newblock   J. Operator Theory \textbf{63} (1) (2010) 47-83.








 \bibitem{C-E-SCH}
J. De Canni\`ere, M. Enock and J.-M. Schwartz,
 \newblock \textit{Sur deux r\'esultats d’analyse harmonique non-commutative: Une application de la th\'eorie des alg\`ebres de Kac},
\newblock J. Operator Theory,  \textbf{5} (1981) 171-194.







\bibitem{EFFROSRUAN1991}
E.G. Effros, Z.-J. Ruan,
\newblock  \textit{A new approach to operator spaces},
\newblock  Canad. J. Math. \textbf{34} (1991) 329-337.



\bibitem{EFFROSRUAN2000}
E.G. Effros, Z.-J. Ruan,
\newblock  \textit{Operator Spaces},
\newblock  London Math. Soc. Monogr. (N.S.) \textbf{23} Clarendon Press, Oxford University Press, New York, 2000.




\bibitem{EYMARD1964}
 P. Eymard,
\newblock \textit{ L'alg\'ebre de Fourier d'un groupe localement compact},
\newblock  Bull. Soc. Math. France \textbf{92} (1964) 181-236.




\bibitem{FIGATALAMANCA1965}
 A. Fig\`a-Talamanca,
\newblock  \textit{Translation invariant operators in $L_p$},
\newblock  Duke Math. J. \textbf{32} (1965) 495-501.






\bibitem{FOLLAND1995}
G. B. Folland,
\newblock \textit{A course in abstract harmonic analysis},
\newblock CRC PRESS, 1995.














\bibitem{GARDELLA2019}
E. Gardella,
\newblock  \textit{A modern look at algebras of operators on $L_p$-spaces},
\newblock   ArXiv:1909.12096v1.




\bibitem{HERZ1974}
C. Herz,
\newblock  \textit{Une g\'en\'eralisation de la notion de transform\'ee de Fourier-Stieltjes},
\newblock   Ann. Inst. Fourier \textbf{24} (1974), 145–157.


\bibitem{HERZ1971}
C. Herz,
\newblock  \textit{The theory of $p$-spaces with an application to convolution operators, its second dual},
\newblock    Trans. Amer. Math. Soc. \textbf{154} (1971) 69-82.









\bibitem{ILIE2013}
M. Ilie,
\newblock \textit{A note on $p$-completely bounded homomorphisms of the Fig\`a-Talamanca-Herz algebras},
\newblock  J. Math. Anal. Appl. \textbf{419} (2014) 273-284




\bibitem{ILIESPRONK2005}
 M. Ilie, N. Spronk,
\newblock  \textit{ Completely bounded homomorphisms of the Fourier algebra}, 
\newblock  J. Funct. Anal. \textbf{225} (2005) 480-499.   
   






   


   
\bibitem{LOSERT1984}
V. Losert,
\newblock  \textit{Properties of the Fourier algebra that are equivalent to amenability},
\newblock  Proc. Amer. Math. Soc. \textbf{92} (1984), 347-354.   





\bibitem{LEMERDY1996}
C. Le Merdy,
\newblock  \textit{Factorization of $p$-completely bounded multilinear maps},
\newblock Paciﬁc J. Math. \textbf{172} (1996), 187-213.












 \bibitem{PAT}
 A. L. T. Paterson,
 \newblock  \textit{The Fourier algebra for locally compact groupoids}, 
 \newblock Canad. J.
Math. \textbf{56} (2004) 1259-1289.



\bibitem{PISIER2001}
G. Pisier,
\newblock  \textit{Similarity Problems and Completely Bounded Maps},
\newblock Second, expanded edition, Lecture Notes in Math., vol. 1618, Springer-Verlag, Berlin 2001.









\bibitem{RUNDE2005}
 V. Runde,
\newblock \textit{Representations of locally compact groups on ${QSL}_p$-spaces and a $p$-analog of the Fourier-Stieltjes algebra},
\newblock 
Pacific J. Math. \textbf{221} (2) (2005) 379-397.







 \bibitem{SHAMS}
 M. Shams Yousefi,
 \newblock \textit{$p$-analog of the semigroup Fourier-Stieltjes algebras},
 \newblock Iranian J. Math. Sci. and Inf.,\textbf{10} (2) (2015), 55-66 


\end{thebibliography}
\end{document}